\documentclass[12pt,savetrees]{amsart}
\usepackage{amssymb, hyperref}

\newcommand{\Bgp}{{\Z^\N}}

\long\def\forget#1\forgotten{}
\newcommand{\issuenumber}{36}
\newcommand{\issuemonth}{November}
\newcommand{\issueyear}{2013}

\newcommand{\sub}{\subseteq}
\newcommand{\R}{\mathbb{R}}
\newcommand{\Cp}{\mathrm{C}_\mathrm{p}}
\newcommand{\Op}{\mathrm{O}}

\newcommand{\ed}{

\section{Unsolved problems from earlier issues}

\begin{issue}Is $\binom{\Omega}{\Gamma}=\binom{\Omega}{\Tau}$?\end{issue}
\begin{issue}Is $\ufin(\Op,\Omega)=\sfin(\Gamma,\Omega)$?And if not, does $\ufin(\Op,\Gamma)$ imply
$\sfin(\Gamma,\Omega)$?\end{issue}
\stepcounter{issue}\begin{issue}Does $\sone(\Omega,\Tau)$ imply $\ufin(\Gamma,\Gamma)$?\end{issue}
\begin{issue}Is $\fp=\fp^*$? (See the definition of $\fp^*$ in that issue.)\end{issue}
\begin{issue}Does there exist (in ZFC) an uncountable set satisfying $\sfin(\cB,\cB)$?\end{issue}
\stepcounter{issue}
\begin{issue}Does $X \nin \NON(\cM)$ and $Y\nin\mathsf{D}$ imply that $X\cup Y\nin \COF(\cM)$?\end{issue}
\begin{issue}[CH]Is $\split(\Lambda,\Lambda)$ preserved under finite unions?\end{issue}
\begin{issue}Is $\cov(\cM)=\fo$? (See the definition of $\fo$ in that issue.)\end{issue}
\stepcounter{issue}
\begin{issue}Could there be a Baire metric space $M$ of weight $\aleph_1$ and a partition
$\mathcal{U}$ of $M$ into $\aleph_1$ meager sets where for each ${\mathcal U}'\subset\mathcal U$,
$\bigcup {\mathcal U}'$ has the Baire property in $M$?\end{issue}
\stepcounter{issue} 
\begin{issue}Does there exist (in ZFC) a set of reals $X$ of cardinality $\fd$ such that all
finite powers of $X$ have Menger's property $\sfin(\Op,\Op)$?\end{issue}
\begin{issue}Can a Borel non-$\sigma$-compact group be generated by a Hurewicz subspace?\end{issue}
\begin{issue}[MA]Is there $X\sbst\bbR$ of cardinality continuum, satisfying $\sone(\BO,\BG)$?\end{issue}
\begin{issue}[CH]Is there a totally imperfect $X$ satisfying $\ufin(\Op,\Gamma)$
that can be mapped continuously onto $\Cantor$?\end{issue}
\begin{issue}[CH]Is there a Hurewicz $X$ such that $X^2$ is Menger but not Hurewicz?\end{issue}
\begin{issue}Does the Pytkeev property of $C_p(X)$ imply that $X$ has Menger's property?\end{issue}
\begin{issue}Does every hereditarily Hurewicz space satisfy $\sone(\BG,\BG)$?\end{issue}
\begin{issue}[CH]Is there a Rothberger-bounded $G\le\Bgp$ such that $G^2$ is not Menger-bounded?\end{issue}
\begin{issue}Let $\cW$ be the van der Waerden ideal. Are $\cW$-ultrafilters closed under products?\end{issue}
\begin{issue}Is the $\delta$-property equivalent to the $\gamma$-property $\binom{\Omega}{\Gamma}$?\end{issue}
\stepcounter{issue}\stepcounter{issue}
\general\end{document}}

\newcommand{\Cantor}{{\{0,1\}^\N}}

\newcommand{\fc}{\mathfrak{c}}\newcommand{\fd}{\mathfrak{d}}
\newcommand{\fp}{\mathfrak{p}}
\newcommand{\NON}{{\mathsf   {NON}}}\newcommand{\COF}{{\mathsf   {COF}}}

\newcommand{\cM}{\mathcal{M}}
\newcommand{\cov}{\mathsf{cov}}

\newcommand{\CH}{the Continuum Hypothesis}\newcommand{\bbR}{\mathbb{R}}
\newcommand{\fo}{\mathfrak{od}}

\renewcommand{\split}{\mathsf{Split}}\newcommand{\bq}{\begin{quote}}\newcommand{\eq}{\end{quote}}
\newcommand{\cB}{\mathcal{B}}\newcommand{\BG}{\cB_\Gamma}
\newcommand{\BO}{\cB_\Omega}

\newcommand{\sone}{\mathsf{S}_1}\newcommand{\sfin}{\mathsf{S}_\mathrm{fin}}
\newcommand{\ufin}{\mathsf{U}_\mathrm{fin}}\newcommand{\gone}{\mathsf{G}_1} 
\newcommand{\nin}{\not\in}
\newcommand{\cW}{\mathcal{W}}

\newcommand{\N}{\mathbb{N}}\newcommand{\Z}{\mathbb{Z}}
\newcommand{\sbst}{\subseteq}
\newcommand{\by}[2]{\par\hfill\emph{#1}, #2}\newcommand{\nby}[1]{\par\hfill\emph{#1}}\newcommand{\Tau}{\mathrm{T}}
\newcommand{\CE}{\textsc{CE}}
\newtheorem{thm}{Theorem}[section]\newcommand{\bthm}{\begin{thm}} \newcommand{\ethm}{\end{thm}}
\newtheorem{prop}[thm]{Proposition}\newcommand{\bprp}{\begin{prop}} \newcommand{\eprp}{\end{prop}}
\newtheorem{fact}[thm]{Fact}\newcommand{\bfct}{\begin{fact}} \newcommand{\efct}{\end{fact}}
\newtheorem{prob}[thm]{Problem}\newcommand{\bprb}{\begin{prob}} \newcommand{\eprb}{\end{prob}}
\newtheorem{lem}[thm]{Lemma}\newcommand{\blem}{\begin{lem}} \newcommand{\elem}{\end{lem}}
\newtheorem{claim}[thm]{Claim}\newcommand{\bclm}{\begin{claim}} \newcommand{\eclm}{\end{claim}}
\newtheorem{cor}[thm]{Corollary}\newcommand{\bcor}{\begin{cor}} \newcommand{\ecor}{\end{cor}}
\newtheorem{conj}[thm]{Conjecture}\newcommand{\bcnj}{\begin{conj}} \newcommand{\ecnj}{\end{conj}}
\theoremstyle{definition}\newtheorem{defn}[thm]{Definition}\newcommand{\bdfn}{\begin{defn}} \newcommand{\edfn}{\end{defn}}
\theoremstyle{remark}\newtheorem{rem}[thm]{Remark}\newcommand{\brem}{\begin{rem}} \newcommand{\erem}{\end{rem}}
\newtheorem{cnv}[thm]{Convention}\newcommand{\bcnv}{\begin{cnv}} \newcommand{\ecnv}{\end{cnv}}
\newtheorem{exam}[thm]{Example}\newcommand{\bexm}{\begin{exam}} \newcommand{\eexm}{\end{exam}}
\newtheorem{issue}{Issue}\newcommand{\bpf}{\begin{proof}} \newcommand{\epf}{\end{proof}}
\newcommand{\be}{\begin{enumerate}}\newcommand{\ee}{\end{enumerate}}\newcommand{\bi}{\begin{itemize}}
\newcommand{\ei}{\end{itemize}}
\newcommand{\general}{\small\vfill\par\noindent\hrulefill\par
\noindent\textbf{Previous issues.} The previous issues of this
bulletin are available online at\\
\url{http://front.math.ucdavis.edu/search?\&t=\%22SPM+Bulletin\%22}
\\[0.1cm]
\textbf{Contributions.} Announcements, discussions, and open problems should be emailed
to \texttt{tsaban@math.biu.ac.il}\\[0.1cm]
\textbf{Subscription.}
To receive this bulletin (free) to your e-mailbox, e-mail us.
}

\newcommand{\arXivl}[4]{\subsection{#2}{#4}\par\hfill{\arx{#1}}\par\hfill\emph{#3}}
\newcommand{\arXiv}[3]{\subsection{#2}\mbox{}\par\hfill{\arx{#1}}\par\hfill\emph{#3}}

\newcommand{\AMS}[3]{\subsection{#1}\mbox{}\par\hfill{\texttt{#3}}\par\hfill\emph{#2}}

\newcommand{\arx}[1]{\url{http://arxiv.org/abs/#1}}

\title[$\mathcal{SPM}$ Bulletin \textbf{\issuenumber} (\issuemonth{} \issueyear)]{%
$\mathcal{SPM}$ Bulletin\\[0.5cm]
Issue number \issuenumber: \issuemonth{} \issueyear{} \CE{}}

\begin{document}
\maketitle


\forget
\section{Editor's note}

\medskip

With best regards,

\by{Boaz Tsaban}{tsaban@math.biu.ac.il}

\hfill \texttt{http://www.cs.biu.ac.il/\~{}tsaban}
\forgotten

\section{Proceedings of the Fourth Workshop on Coverings, Selections and Games in Topology}

Volume 160, Issue 18, of Topology and its Applications (1 December 2013, Pages 2233--2566), is dedicated
to the proceedings of the Fourth Workshop on Coverings, Selections and Games in Topology, 
on the occasion of Ljubi\v{s}a D.R. Ko\v{c}inac's 65$^\mathrm{th}$ birthday.
This issue is now available online at

\centerline{\url{http://www.sciencedirect.com/science/journal/01668641/160/18}}

Contents:
\be
\item \emph{Preface}, Page 2233, G. Di Maio, B. Tsaban.

\item \emph{The mathematics of Ljubi\v{s}a D.R. Ko\v{c}inac},
Pages 2234--2242,
B. Tsaban.

\item \emph{Selectively c.c.c.\ spaces}, Pages 2243-2250, L.F. Aurichi.

\item \emph{Weak covering properties and selection principles},
Pages 2251--2271, L. Bab\-inkostova, B.A. Pansera, M. Scheepers.
 
\item \emph{Detecting topological groups which are (locally) homeomorphic to LF-spaces},
Pages 2272--2284, T. Banakh, K. Mine, D. Repov\v{s}, K. Sakai, T. Yagasaki.
 
\item \emph{Gap topologies in metric spaces}, Pages 2285-2308, G. Beer, C. Costantini, S. Levi.
 
\item \emph{On two selection principles and the corresponding games}, Pages 2309--2313,
A. Bella.
 
\item \emph{Scale function vs topological entropy},  Pages 2314--2334, 
F. Berlai, D. Dikranjan, A. Giordano Bruno.

\item \emph{On diagonal resolvable spheres},
Pages 2335--2339, 
P.V.M. Blagojevi\'c, A.D. Blagojevi\'c, L.D.R. Ko\v{c}inac.
 
\item \emph{The Schwarz genus of the Stiefel manifold},
Pages 2340--2350, P.V.M. Blagojevi\'c, R.N. Karasev.
 
\item \emph{Dowker-type example and Arhangelskiiʼs image-property},
Pages 2351--2355,
M. Bonanzinga, M. Matveev.
 
\item \emph{Modifications of sequence selection principles},
Pages 2356--2370,
L. Bukovsk\'y, J. \v{S}upina.
 
\item \emph{On the cardinality of the $\theta$-closed hull of sets},
Pages 2371--2378,
F. Cammaroto, A. Catalioto, B.A. Pansera, B. Tsaban.
 
\item \emph{Variations on selective separability in non-regular spaces},
Pages 2379--2385,
A. Caserta, G. Di Maio.
 
\item \emph{Some generalizations of Backʼs Theorem},
Pages 2386--2395,
A. Caterino, R. Ceppitelli, L. Hol\'a.
 
\item \emph{On Haar meager sets},
Pages 2396--2400,
U.B. Darji.
 
\item \emph{Some further results on image-γ and image-image-covers},
Pages 2401--2410,
P. Das, D. Chandra.
 
\item \emph{Indestructibility of compact spaces},
Pages 2411--2426,
R.R. Dias, F.D. Tall.
 
\item \emph{On characterized subgroups of compact abelian groups},
Pages 2427--2442,
D. Dikranjan, S.S. Gabriyelyan.
 
\item \emph{Productively Lindelöf and indestructibly Lindelöf spaces},
Pages 2443--2453,
H. Duanmu, F.D. Tall, L. Zdomskyy.
 
\item \emph{Universal frames},
Pages 2454--2464,
T. Dube, S. Iliadis, J. van Mill, I. Naidoo.
 
\item \emph{Continuous weak selections for products},
Pages 2465--2472,
S. Garc\'{\i}a--Ferreira, K. Miyazaki, T. Nogura.
 
\item \emph{Selections and metrisability of manifolds},
Pages 2473--2481,
D. Gauld.
 
\item \emph{On base dimension-like functions of the type Ind},
Pages 2482--2494,
D.N. Georgiou, S.D. Iliadis, A.C. Megaritis.
 
\item \emph{Selection properties of uniform and related structures},
Pages 2495--2504,
L.D.R. Ko\v{c}inac, H.P.A. K\"unzi.
 
\item \emph{A characterization of the Menger property by means of ultrafilter convergence},
Pages 2505--2513,
P. Lipparini.
 
\item \emph{Transfinite extension of dimension function image},
Pages 2514--2522,
N.N. Martynchuk.
 
\item \emph{Homogeneity and h-homogeneity},
Pages 2523--2530,
S.V. Medvedev.
 
\item \emph{The weak Hurewicz property of Pixley–Roy hyperspaces},
Pages 2531--2537,
M. Sakai.
 
\item \emph{Comparing weak versions of separability},
Pages 2538--2566,
D.T. Soukup, L. Soukup, S. Spadaro.
\ee

We use this opportunity to thank the Chief Editors and the referees for their tremendous help.

\nby{Giuseppe Di Maio, Boaz Tsaban}

\hfill{Guest Editors, Topology and its Applications}

\section{Long announcements}

\arXivl{1304.6628}
{Completeness and related properties of the graph topology on function spaces}
{Lubica Hol\'a, L\'aszl\'o Zsilinszky}
{The graph topology $\tau_{\Gamma}$ is the topology on the space $C(X)$ of all
continuous functions defined on a Tychonoff space $X$ inherited from the
Vietoris topology on $X\times \mathbb R$ after identifying continuous functions
with their graphs. It is shown that all completeness properties between
complete metrizability and hereditary Baireness coincide for the graph topology
if and only if $X$ is countably compact; however, the graph topology is
$\alpha$-favorable in the strong Choquet game, regardless of $X$. Analogous
results are obtained for the fine topology on $C(X)$. Pseudocompleteness, along
with properties related to 1st and 2nd countability of $(C(X),\tau_{\Gamma})$
are also investigated.}

\arXivl{1306.5463}
{Topological games and Alster spaces}
{Leandro F. Aurichi and Rodrigo R. Dias}
{In this paper we study connections between topological games such as
Rothberger, Menger and compact-open, and relate these games to properties
involving covers by $G_{\delta}$ subsets. The results include: (1) If Two has a
winning strategy in the Menger game on a regular space $X$, then $X$ is an Alster
space. (2) If Two has a winning strategy in the Rothberger game on a
topological space $X$, then the $G_\delta$-topology on $X$ is Lindelof. (3) The
Menger game and the compact-open game are (consistently) not dual.}

\arXivl{1306.5421}
{Strongly Summable Ultrafilters, Union Ultrafilters, and the Trivial Sums
  Property}
{David J. Fern\'andez Bret\'on}
{We answer two questions of Hindman, Stepr\=ans and Strauss, namely we prove
that every strongly summable ultrafilter on an abelian group is sparse and has
the trivial sums property. Moreover we show that in most cases the sparseness
of the given ultrafilter is a consequence of its being isomorphic to a union
ultrafilter. However, this does not happen in all cases: we also construct
(assuming $\cov(\cM)=\fc$), on the Boolean group, a strongly summable ultrafilter that
is not additively isomorphic to any union ultrafilter.}

\arXivl{1307.7928}{Topological games and productively countably tight spaces}
{Leandro F. Aurichi and Angelo Bella}
{The two main results of this work are the following: if a space $X$ is such
that player II has a winning strategy in the game $\gone(\Omega_x, \Omega_x)$
for every $x \in X$, then $X$ is productively countably tight. On the other
hand, if a space is productively countably tight, then $\sone(\Omega_x,
\Omega_x)$ holds for every $x \in X$. With these results, several other results
follow, using some characterizations made by Uspenskii and Scheepers.}

\arXivl{1310.0409}
{On a game theoretic cardinality bound}
{Leandro F. Aurichi and Angelo Bella}
{The main purpose of the paper is the proof of a cardinal inequality for a
space with points $G_\delta$, obtained with the help of a long version of the
Menger game. This result improves a similar one of Scheepers and Tall.}

\arXivl{1311.1468}
{Pytkeev $\aleph_0$-spaces}
{Taras Banakh}
{A regular topological space $X$ is defined to be a Pytkeev $\aleph_0$-space
if it has countable Pytkeev network. A family $\mathcal P$ of subsets of a
topological space $X$ is called a Pytkeev network in $X$ if for a subset
$A\subset X$, a point $x\in{\bar A}$ and a neighborhood $O_x\subset X$ there is
a set $P\in\mathcal P$ such that $x\in P\subset O_x$ and moreover $P\cap A$ is
infinite if $x$ is an accumulation point of $A$. The class of Pytkeev
$\aleph_0$-spaces contains all metrizable separable spaces and is (properly)
contained in the class of $\aleph_0$-spaces. This class is closed under many
operations over topological spaces: taking a subspace, countable Tychonoff
product, small countable box-product, countable direct limit, the hyperspace.
For an $\aleph_0$-space $X$ and a Pytkeev $\aleph_0$-space $Y$ the function
space $C_k(X,Y)$ endowed with the compact-open topology is a Pytkeev
$\aleph_0$-space. A topological space is second countable if and only if it is
Pytkeev $\aleph_0$-space with countable fan tightness.}

\arXivl{1310.8622}
{Selective covering properties of product spaces, II: $\gamma$ spaces}
{Arnold W. Miller, Boaz Tsaban, Lyubomyr Zdomskyy}
{We study productive properties of  $\gamma$ spaces, and their relation to other, classic and modern, selective covering properties.
Among other things, we prove the following results:
\be
\item Solving a problem of F. Jordan, we show that for every unbounded tower set $X\sub\R$ of cardinality $\aleph_1$,
the space $\Cp(X)$ is productively Fr\'echet--Urysohn. In particular, the set $X$ is productively $\gamma$.
\item Solving problems of Scheepers and Weiss, and proving a conjecture of Bab\-inkostova--Scheepers, we prove that,
assuming \CH{}, there are  $\gamma$ spaces whose product is not even Menger.
\item Solving a problem of Scheepers--Tall, we show that
the properties $\gamma$ and Gerlits--Nagy (*) are preserved by Cohen forcing. Moreover, every Hurewicz space that Remains
Hurewicz in a Cohen extension must be Rothberger (and thus (*)).
\ee
We apply our results to solve a large number of additional problems,
and use Arhangel'\-ski\u{\i} duality to obtain results concerning local properties of function spaces and countable topological
groups.}

\arXivl{1311.2011}
{Productively countably tight spaces of the form $C_k(X)$}
{Leandro Fiorini Aurichi, Renan Maneli Mezabarba}
{Some results in $C_k$-theory are obtained with the use of bornologies. We
investigate under which conditions the space of the continuous real functions
with the compact-open topology is a productively countably tight space, which
yields some applications on Alster spaces.}

\section{Short announcements}\label{RA}

\arXiv{1304.0472}
{Partitioning bases of topological spaces}
{Daniel T. Soukup, Lajos Soukup}

\arXiv{1304.0486}
{Compactness of $\omega^\lambda$}
{Paolo Lipparini}

\arXiv{1304.2042}
{Non-meager free sets for meager relations on Polish spaces}
{Taras Banakh and Lyubomyr Zdomskyy}

\AMS{On strong $P$-points}
{Andreas Blass; Michael Hrusak; Jonathan Verner}
{http://www.ams.org/journal-getitem?pii=S0002-9939-2013-11518-2}

\arXiv{1307.0184}
{Products and countable dense homogeneity}
{Andrea Medini}

\arXiv{1307.1989}
{Strong colorings yield kappa-bounded spaces with discretely untouchable
  points}
{Istvan Juhasz and Saharon Shelah}

\arXiv{1310.1827}{Selective and Ramsey ultrafilters on $G$-spaces}
{O.V. Petrenko, I.V. Protasov}

\arXiv{1311.1677}
{Seven characterizations of non-meager P-filters}
{Kenneth Kunen, Andrea Medini, Lyubomyr Zdomskyy}

\ed